\documentclass{amsart}
\usepackage{amsfonts,amscd}

\newtheorem{theorem}{Theorem}[section]
\newtheorem{lemma}[theorem]{Lemma}
\newtheorem{corollary}[theorem]{Corollary}
\newtheorem{proposition}[theorem]{Proposition}
\theoremstyle{remark}

\theoremstyle{definition}

\numberwithin{equation}{section}
\makeatother

\DeclareMathOperator{\Cdb}{{\mathbb C}}

\DeclareMathOperator{\Tdb}{{\mathbb T}}

\begin{document}

\title[Duality and operator algebras II]{Duality and operator algebras II: \\
 Operator algebras as Banach 
algebras}
 
\date{August 20, 2004.  Revised October 21, 2004.}
 
\author{David P. Blecher}
\address{Department of Mathematics, University of Houston, Houston,
TX
77204-3008}
\email[David P. Blecher]{dblecher@math.uh.edu}
 \author{Bojan Magajna} 
\address{Department of Mathematics,
University of Ljubljana,
Jadranska 19, Ljubljana 1000, Slovenia}
\email[Bojan Magajna]{Bojan.Magajna@fmf.uni-lj.si}
\thanks{*Blecher is partially supported by a grant from
the National Science Foundation.
Magajna is  partially supported by the Ministry of Science and
Education of Slovenia. } 

\begin{abstract}    We answer, by counterexample,
several  open questions concerning 
algebras of operators on a Hilbert
space.  The answers add further weight to the thesis that, for many 
purposes, such algebras 
ought to be studied in the framework of operator spaces,
as opposed to that of Banach spaces and Banach algebras.
In particular, the `nonselfadjoint analogue' of a W*-algebra
resides naturally in the category of dual operator
spaces, as opposed to dual Banach spaces.
We also show that an automatic w*-continuity result  
in the preceding paper of the authors, is sharp.
\end{abstract}
 
\maketitle

\section{Introduction}

An {\em operator algebra} is an algebra of operators on a
Hilbert space.    
Since the advent of operator space theory, there has been much 
progress toward the development of a {\em general} theory of 
such algebras (e.g.\ see \cite{BLM}).   From the `operator space perspective',
on an operator algebra $A$ one should 
consider not only the norm on $A$, but also the 
canonical norms on the spaces $M_n(A)$ of matrices with entries in $A$,
for all $n \geq 1$.    The obvious question is, is this really 
necessary for the study of such algebras? 
  While admittedly this is not a well-posed question,
since it depends on the applications one has in mind,
a recent survey \cite{Bare} collected some  
test questions which have resisted solution to date,
and  whose answers might 
`tip the balance' on this issue, in some sense.  We are now able to answer 
several of these questions.   Our main result may be summarized as 
saying that 
Sakai's famous characterization
of von Neumann algebras 
in terms of $C^*$-algebras with a Banach space predual
(see \cite[Theorem 1.16.7]{Sakai}), is not valid for 
general operator algebras
without using the operator space framework. In particular, we exhibit here an
operator algebra with an identity of norm 1,  even a 
commutative one, which has a Banach space predual, but is not   
homomorphic, via a homeomorphism for the associated               
weak* topologies, to any $\sigma$-weakly closed (that is, weak* closed) 
 operator algebra.   This rules out the possibility, which 
had remained open,
of a `non-operator-space variant' of the following theorem attributable
to Le Merdy and the two authors (e.g.\ see \cite[Section 2.7]{BLM} and 
\cite{LM,BM}): namely that
the  $\sigma$-weakly closed operator  algebras `are precisely' the 
operator algebras which possess an operator space predual.
Thus, we are able to bring to its final form
the topic of abstract characterizations
of $\sigma$-weakly closed 
operator algebras.
We also use our counterexample to deduce  
that several other 
known results about operator algebras and 
operator spaces are not valid if one drops hypotheses
involving `matrix norms'.      For example, 
we exhibit a subspace $X$ of a $C^*$-algebra $A$, and an
$a \in A$ with  $a X \subset X$, such that
$X$ is isometric to a dual
Banach space, but the function $x \mapsto a x$ on $X$
is not weak$^*$ continuous (we showed in \cite{BM} that this 
function is always weak$^*$ continuous if also
$a^* X \subset X$, or if $X$ is a dual operator space).  
As another byproduct, we have found  more simple examples of 
operator spaces which have a Banach space predual
which is not an operator space predual (see \cite{BM} for more discussion
of this point).     

We refer the reader to any of the recently available texts on
operator spaces, for more information on that topic if needed. 
For the duality of operator 
spaces, we recommend \cite[Section 1.4]{BLM}, although this is
not necessary for reading our paper.  
We abbreviate `weak*' to `w*-' throughout.

\section{Duality and lowersemicontinuity}
To construct 
`noncanonical' Banach space preduals, the following well known result
is very useful.  A function defined on a  
dual Banach space will be called
{\em lower w*-continuous} if it is lowersemicontinuous with respect to the
w*-topology.

\begin{lemma} \label{ERlem}  {\rm (\cite[Lemma 3.1]{ER}, \cite[Lemma 8.8]{F})}  
If $X$ is a dual Banach space, and if $| | | \cdot | | |$
 is an equivalent norm on $X$, then the following
are equivalent:
\begin{itemize} \item [(i)]  $\{ x \in X : | | | x 
| | | \leq 1 \}$ is w*-closed;
  \item [(ii)]   $| | | \cdot | | |$ 
is lower w*-continuous;
\item [(iii)]  $(X, | | | \cdot | | |)$ is a dual Banach space,
and the associated w*-topology is the same
as the original w*-topology on $X$.   
\end{itemize}
\end{lemma}

Let $Y$ be an operator space which is also a
non-reflexive dual Banach space, with predual $Y_*$.  
In fact, in our examples, $Y$ will 
be a subspace of a unital $C^*$-algebra $A$ with $1_A \in Y$. 
We will identify $\Cdb$ with $\Cdb 1_A \subset A$.
Suppose that $T$ is a bounded operator on $Y$,
which is discontinuous with respect to the w*-topology of $Y$.
We then consider operator spaces which may be `built' 
from $Y$ and $T$.  In particular, in the remainder of the 
paper we will be using the 
the following four subspaces of $M_2(A)$:
$$ B \; = \; \Bigl\{ \left[\begin{array}{ccl} 0 & y\\ 0 & Ty\end{array}\right]
\, , \, y \in Y \Bigr\} , $$ 
$$ C  \; = \; \Bigl\{ \left[\begin{array}{ccl} a & y\\ 0 & a + Ty\end{array}\right]
\, , \, y \in Y , a \in \Cdb 1_A \Bigr\} , $$
$$E  \; = \; \Bigl\{
\left[\begin{array}{ccl} Ty & y\\ 0 & Ty\end{array}\right]
\, , \, y \in Y \Bigr\} , $$
and 
$$D \; = \; \Bigl\{ \left[\begin{array}{ccl} a + Ty & y \\ 0 & a + Ty
\end{array}\right] \, , \, y \in Y , a \in \Cdb 1_A  \Bigr\} . $$
Notice that $B, C,$ and $D$ are  operator algebras, 
subalgebras of $M_2(A)$, if $T$ takes 
values in $\Cdb 1_A$.   
The norms of the matrices in the  expressions for $B$ and $E$ above,
are 
equivalent to the original norm on $Y$, whereas the norms
in the  expressions for $C$ and $D$ are equivalent to the $\infty$-direct 
sum norm on
$Y \oplus \Cdb$.  Our strategy will be as follows.
Consider for example the canonical isomorphism
$\theta : Y \oplus^\infty \Cdb \to D$ given by
$$\theta((y,a))  \; = \; \left[\begin{array}{ccl} a + Ty & y \\ 0 & a + Ty
\end{array}\right] \; , \qquad y \in Y , a \in \Cdb .$$
The $M_2(A)$-norm on $D$ transfers, via $\theta$, to a
norm $| | | \cdot | | |$ on the dual space $Y \oplus^\infty \Cdb$.
By Lemma \ref{ERlem}, $D$ is a dual Banach space if 
$| | | \cdot | | | = 
\Vert \theta(\cdot) \Vert_{D}$ is lower w*-continuous with respect to 
the w*-topology induced by the predual $Y_* \oplus^1 \Cdb$.
Similar (but simpler)
considerations apply to the other spaces $B, C, E$ above.
Thus $B$ is a dual Banach space and the canonical map
$\rho : Y \to B$ is a w*-homeomorphism,
 if $\Vert \rho(\cdot) \Vert_{B}$ is 
lower w*-continuous on $Y$.

For computing the norms above, we use 
the following explicit formula:
\begin{equation} \label{eqtr}
\Bigl\Vert \left[\begin{array}{ccl} a & x \\ 0 & b \end{array}\right]
\Bigr\Vert^2 = \frac{1}{2} [|a|^2 + \Vert x \Vert^2 + |b|^2 
+ \sqrt{(|a|^2 + \Vert x \Vert^2 + |b|^2)^2 - 4 |a|^2 |b|^2} ] , 
\end{equation} 
for $x \in A$ and $a, b \in \Cdb 1_A$.  This follows easily 
from the last formula in the proof of 2.2.11 in
\cite{BLM}.  
 
Because the computations are quite manageable here,
we will henceforth restrict our attention to the case when 
$Y = \ell^1$ and $Y_* = c_0$.   We write $(e_k)$ for the canonical
basis.  We will fix an operator space structure 
on $Y$ so that the identity of the containing $C^*$-algebra $A$ is 
$e_1$, which we will write as $1_Y$.  For example, 
let $A$ be the commutative $C^*$-algebra $C(\Tdb^\infty)$, where 
$\Tdb$ is the unit circle, and $Y$ the copy of $\ell^1$ in $A$
corresponding to the functions $(z_k) \mapsto a_1 + \sum_{k\geq 2} a_k z_k$
on $\Tdb^\infty$, for all $(a_k) \in \ell^1$.
(See also e.g.\ \cite[9.6]{Pisier} or \cite[4.3.8]{BLM}.) 
Note that 
w*-convergence of bounded nets in $\ell^1$
simply means `component-wise convergence'.
We will take the map $T$ above to be of the form $Ty=\tau(y)1$, where $\tau\in      
Y^*\setminus Y_*$ is of norm $1$.
 We will need the following:  

\begin{lemma}\label{l1} For a functional $\tau=(a_j)\in\ell^{\infty}=(\ell^1)^*$,
 the function $f(y)=\|y \|^2+|\tau(y)|^2$ is lower w*-continuous on $\ell^{1}$
(regarded as the dual of $c_0$) if and only if
\begin{equation}\label{MS} MS\leq1,\ \mbox{where}\ M=\sup_j|a_j|\ \mbox{and}\ S=\limsup_j|a_j|.
\end{equation}
\end{lemma}

\begin{proof}  If $f$ is lower w*-continuous, consider the sequence
$y(m)=ze_1+we_m$, where 
$z,w\in \Cdb$.
This sequence converges in the w*-topology 
to $y=ze_1$.  Lower w*-continuity then demands that 
$f(y)\leq\liminf_mf(y(m))$.  This can be rewritten as
\begin{equation}\label{21} 0\leq\liminf_m[(|a_m|^2+1)|w|^2+2(|z||w|+{\rm Re}(a_1\overline{a}_m
z\overline{w}))].
\end{equation}
If $a$ is any limit point of the sequence $(a_m)$, then by
choosing first a subsequence converging
to $a$, and then appropriate phases of $z$ and $w$, it follows from (\ref{21}) that
$0\leq(|a|^2+1)t^2+2st(1-|a_1||a|)$ for all $s,t>0$.
Letting $t\to0$, we have $1-|a_1||a|\geq0$.
In particular $|a_1|S\leq 1$. Similarly, $|a_k|S\leq1$ for all $k$.
Hence $MS\leq1$. 

Conversely, note that $f(y)=|\sum_{j=1}^{\infty}a_jy_j|^2+(\sum_{j=1}^{\infty}|y_j|)^2$ for $y=(y_j)\in\ell^1$, so 
\begin{equation}\label{22} f(y)=
\sum_{j=1}^{\infty}(1+|a_j|^2)|y_j|^2+2\sum_{i<j}(|y_i||y_j|+{\rm Re}(a_i\overline{a}_jy_i
\overline{y}_j)) .
 \end{equation}
Of course, 
$$ |y_i||y_j|+{\rm Re}(a_i\overline{a}_jy_i\overline{y}_j)\geq
|y_i||y_j|(1-|a_i||a_j|).
$$
We claim that if $MS<1$, then $|a_i||a_j|>1$ only for
finitely many pairs $(i,j)$, so that almost all the terms on the right side of 
(\ref{22}) are non-negative. Indeed, if this were not true, then
 there are two possibilities:
(1) there exists an $i$ such that $|a_i||a_j|>1$ for infinitely many $j$'s,
 or (2) there
exist infinitely many $i$'s such that for some $j(i)$ we have $|a_i||a_{j(i)}|
>1$.  In the
first case, it follows that $|a_i|S\geq1$, hence $MS\geq1$, which contradicts the assumption. 
Similarly, in the second case we have $|a_i|M>1$, hence $SM\geq1$, again a
contradiction. 

Thus if $MS<1$, the sum
 $\sum_{i<j} |y_i||y_j|+{\rm Re}(a_i\overline{a}_jy_i
\overline{y}_j)$ splits into a finite partial sum, 
and an infinite sum in which 
all the terms $|y_i||y_j|+{\rm Re}(a_i\overline{a}_j y_i\overline{y}_j)$
are nonnegative.  The finite partial sum is actually $w^*$-continuous.
The other sum is the supremum of its own finite partial sums, each of 
which is $w^*$-continuous.  Since a supremum of lowersemicontinuous
functions is lowersemicontinuous, this proves that 
our infinite sum is lower $w^*$-continuous.  A similar but 
easier argument shows that $\sum_{j=1}^{\infty}(1+|a_j|^2)|y_j|^2$
is lower $w^*$-continuous.  

If $MS=1$, we note that $f$ is the supremum of functions $f_n(y)=
\|y\|^2+|\tau_n(y)|^2$, where $\tau_n=(1-\frac{1}{n})\tau$ ($n=1,2,\ldots$). 
Since the functions $f_n$ are
lower $w*$-continuous by what we have already proved, so must be $f$. 
\end{proof}

\section{Consequences}

\begin{corollary} \label{idemp}
There exists an operator algebra $B$, which is a dual Banach space,
and an idempotent element
$p\in B$, such that ${\rm dim}(pB) = 1$, and
such that left multiplication by $p$ is not w*-continuous on $B$.
\end{corollary}

\begin{proof}  We use the notation in the previous section.
Choose $\tau$ to satisfy the conditions of Lemma \ref{l1}, and $\tau(1_Y) = a_1 = 1$.
Let $T = \tau(\cdot) 1_Y$, and let $B$ be as in the 
last section, with $\rho : Y \to B$ the canonical map. 
The  $M_2(A)$-norm on $B$ is given
by the expression $\Vert \rho(y) \Vert_B 
= \sqrt{\Vert y \Vert^2 + |\tau(y)|^2}$ for 
$y \in Y$ (by e.g.\ (\ref{eqtr})).
By Lemma \ref{l1}, this quantity is lower w*-continuous.  It follows,
as in the arguments a couple of  paragraphs
above Lemma \ref{l1}, that $B$ is a 
dual Banach space, with w*-topology determined by the canonical
pairing with $c_0$.  If $p \in B$ corresponds to $y = e_1 =
1_Y \in Y = \ell^1$,
then $p$ is idempotent, $p B$ is one-dimensional, and 
the map $b \to pb$ on $B$ corresponds to the map $T$ on $Y$,
which is not w*-continuous.
\end{proof}  

Such an algebra $B$ is not isomorphic (in the appropriate sense)
to a $\sigma$-weakly closed operator algebra,
since
the product on any $\sigma$-weakly closed operator algebra
is separately w*-continuous. 
 To obtain a `unital' counterexample is a little harder:
  
\begin{theorem} \label{main}
There exists a commutative operator algebra $D$ with an identity 
of norm $1$, which is a dual Banach space,
and a nilpotent element $p\in D$, such that left multiplication by $p$ is not
w*-continuous on $D$.  Moreover, $D$ is not 
homomorphic, via a homeomorphism for the associated
weak* topologies, to any $\sigma$-weakly closed operator algebra.
\end{theorem}

\begin{proof}   
The last assertion here follows as in the line
above the theorem.

We employ a similar strategy to that of Corollary \ref{idemp},
and the notation in the previous section.  
Set $\tau_0=(0,1,1,1,\ldots)\in\ell^{\infty}$,
$\tau=\frac{1}{4}\tau_0$, and $T = \tau(\cdot) 1_Y$.
Consider the unital operator algebra $D$ 
above, which in this case consists of all matrices 
\begin{equation}\label{23}x=\left[\begin{array}{ll}
b&y\\
0&b
\end{array}\right],\ \mbox{where}\ y\in Y \ \mbox{and}\ b\in
\Cdb.
\end{equation}
Clearly $D$ contains the space $E$ as a subspace.
Let $p \in E$ correspond to $y = e_1 = 1_Y \in Y = \ell^1$.
We claim that it suffices to prove that the unit ball for the norm
$| | | \cdot | | |$ 
on $Y \oplus \Cdb$,
defined a couple of  paragraphs
above Lemma \ref{l1}, is w*-closed in the topology
given by the 
canonical pairing with $c_0 \oplus \Cdb$. 
(This corresponds to the pairing of 
a matrix $x$ of the form (\ref{23}) and an element
$v=(z,\beta)\in c_0\oplus\mathbb{C}$, via the formula  
$\langle x,v\rangle=\langle y,z\rangle_{\ell^1,c_0}+
\beta(b-\tau(y))$.)
If this is the case, then $D$ is a dual Banach space by 
Lemma \ref{ERlem}, 
and $E$ is a w*-closed subspace, but left multiplication by $p$ is not
w*-continuous on $E$, and hence not on $D$, since $\tau\notin c_0$.

By (\ref{eqtr}), the norm of a matrix of the form (\ref{23}) is
 $\|x\|^2=\frac{1}{2}[2|b|^2+\|y\|^2+\sqrt{\|y\|^4+4|b|^2\|y\|^2}]$. From this
(or otherwise) it follows by elementary algebraic manipulations that
\begin{equation}\label{24} \|x\|\leq1 \; \Longleftrightarrow 
 \; \|y\|+|b|^2\leq1.
\end{equation}
In the topology described a couple of paragraphs
above Lemma \ref{l1}, the convergence of nets in $D$
 is as follows:
$$\left[\begin{array}{ccl} a_t + Ty_t & y_t \\ 0 & a_t + T y_t
 \end{array}\right]  \to \left[\begin{array}{ccl} a + Ty & y \\ 0 & a + Ty
\end{array}\right]  \iff a_t \to a \; \textrm{in} \Cdb ,
y_t  \to y \, w^* \, \textrm{in} \,  Y .
$$
To prove that the unit ball of $D$
is closed in this topology,
by (\ref{24}) it suffices to show that the function
$$\left[\begin{array}{ccl} a + Ty & y \\ 0 & a + Ty
\end{array}\right]  \; \longmapsto  \; \|y\|+|a+\tau(y)|^2$$
is lower w*-continuous on Ball$(D)$, in the latter topology.
Since Ball$(D)$ is bounded and $a\in\Cdb$, the proof reduces
to showing that for a fixed $a\in\Cdb$, the function
$$g(y)=\|y\|+|a+\tau(y)|^2=|a|^2+\|y\|+\frac{1}{2}{\rm Re}(\overline{a}\tau_0(y))+
\frac{1}{16}|\tau_0(y)|^2$$
is lower w*-continuous on $Y$. We can rewrite $g(y)$ in the form
$$|a|^2+\left(\frac{7}{8}\|y\|+\frac{1}{2}{\rm Re}(\overline{a}
\tau_0(y))\right)+\frac{1}{4}\left(\frac{1}{2}\|y\|-\frac{1}{4}\|y\|^2\right)+\frac{1}{16}
\left(\|y\|^2+|\tau_0(y)|^2\right). $$
The last term in this expression is lower w*-continuous by Lemma \ref{l1}. Since
the norm is lower w*-continuous and the function $t\mapsto t/2-t^2/4$ is increasing
on the interval $[0,1]$, it follows that the term $\|y\|/2-\|y\|^2/4$ is also lower
w*-continuous. It remains to prove that the second term in the last centered equation
is lower
w*-continuous. But this term can be written as $1/2$ of
$$\frac{7}{4}\|y\|+{\rm Re}(\overline{a}\tau_0(y))=\frac{7}{4}|y_1|+
\sum_{j=2}^{\infty}\left(\frac{7}{4}|y_j|+{\rm Re}(\overline{a}y_j)\right).$$
In the last sum all the terms are nonnegative,
 since $\frac{7}{4}|y_j|+{\rm Re}(\overline{a}y_j)
\geq\frac{7}{4}|y_j|-|y_j|\geq0$.
Hence the lower w*-continuity follows, as in the proof
of Lemma \ref{l1}, by considering first the finite partial sums.
\end{proof}

The reader may wonder why, in the last theorem, we did not use the algebra $C$
from the last section, constructed from the simpler algebra $B$ used 
in the proof of Corollary \ref{idemp}, by simply adjoining the identity 
of $M_2(A)$.  In fact, this construction does not produce a dual 
Banach space.  We may use this observation to show that another 
result which is valid for operator algebras which possess an
operator space predual, fails if we assume only
a Banach space predual:

\begin{corollary}  \label{unitiz}  There exists an 
operator algebra $A$ which possesses a Banach space predual,
such that the unitization  $A^1$
of $A$ (see \cite[Section 2.1]{BLM}), 
possesses no Banach space predual for which the 
embedding of $A$ in $A^1$ is w*-continuous.    In contrast,
 there can exist no such example which possesses an operator space 
predual.    
\end{corollary}       

\begin{proof}  The last assertion follows from 
the characterization of dual operator algebras mentioned in the 
first paragraph of our paper, together with \cite[Proposition 2.7.4]{BLM}.   
For the first assertion, construct $C$ as mentioned a couple of
paragraphs above, using the functional
$\tau = (1,1,\cdots) \in \ell^\infty = (\ell^1)^*$.  
Suppose that $C$ possessed a 
Banach space predual for which the canonical
 embedding of $B$ in $C$ was w*-continuous.   By a variant
of the Krein-Smulian theorem, $B$ would then be 
w*-closed in $C$, and the embedding a w*-homeomorphism.  Let 
$\chi : C \to \Cdb$ be evaluation at the 1-1 entry, 
a contractive homomorphism.  Since the kernel $B$ of $\chi$
is w*-closed, $\chi$ is w*-continuous.  It follows that
the map $c \mapsto (c - \chi(c) 1_C , \chi(c))$ is a 
w*-homeomorphism  from $C$ onto $B \oplus^\infty \Cdb$.
Thus a net $(b_t + \lambda_t 1_C)$ in $C$ w*-converges to
$b + \lambda 1_C$ if and only if $\lambda_t \to \lambda$
in $\Cdb$ and $b_t \to b$ in the 
w*-topology in $B$.  The latter condition simply says that
$y_t \to y$ in the
w*-topology of $\ell^1$, if $y_t$ and $y$ are the corresponding 
(via $\rho$) elements in $Y = \ell^1$.
By Lemma \ref{ERlem}, any closed ball
centered at the origin in $C$ is  w*-closed 
with respect to the topology just described 
on $C$.  However, if $y_m = e_1 - \frac{1}{2} e_m$, and 
$b_m$ is the corresponding element of $B$, then 
$(b_m + 1_C)$ is a net in $C$ with w*-limit $b + 1_C$,
where $b \in B$  corresponds to $e_1 \in Y$.     
But  one may easily check using (\ref{eqtr}) that 
$b + 1_C$ lies outside a ball
centered at the origin in $C$ which contains all the terms
$b_m + 1_C$. 
 \end{proof}

Theorem \ref{main} shows that the following is the best result
along the lines above, which one can
hope for in the `Banach algebra category'.

\begin{proposition} \label{LM}  An operator algebra which has
a Banach space predual, and whose product is separately 
w*-continuous, is isometric via a homomorphism which is also 
 a w*-homeomorphism, to a
 $\sigma$-weakly closed operator algebra.
\end{proposition}

\begin{proof}  This is a remark in the Notes to Section 2.7 in
\cite{BLM}.  Indeed, it follows from Le Merdy's proof in \cite{LM}.     
\end{proof}

We also have the following positive result in the case of a 
Banach space predual:

\begin{proposition} \label{delt} Let $B$ be a unital
operator algebra which is
a dual Banach space.  Then 
$\Delta(B) = B \cap B^*$ is a $W^*$-algebra,
and if $b \in \Delta(B)$ then the maps $a \mapsto ab$ and
$a \mapsto ba$ on $B$, are $w^*$-continuous.
\end{proposition}
 
\begin{proof}    
Let $A = B^{**}$ and let $q : A \to B$
be the canonical projection (the adjoint of the inclusion
$B_* \subset B^*$).   Since $q(1) = 1$,
$q$ takes Hermitian elements to Hermitian elements.
That is, $q$ induces a $w^*$-continuous projection $Q$ of
$\Delta(A)$ onto $\Delta(B)$.    Thus
$\Delta(B)$ is isometric to the dual space
$\Delta(A)/{\rm Ker}(Q)$, and so
$\Delta(B)$ is a $W^*$-algebra.
The last part follows from e.g.\ \cite[Theorem 3.3]{BM}.   
  \end{proof}

Theorem \ref{main} also yields solutions
to a couple of other interesting  questions,
as we discuss next.
  
A famous theorem of Tomiyama characterizes `conditional 
expectations' on $C^*$-algebras as the contractive projections
(e.g.\ see \cite[Remark 2.6.5]{Sakai}).  There is 
a known analogue of this for 
nonselfadjoint operator algebras, but it applies to 
{\em completely contractive} projections (see \cite[Corollary 4.2.9]{BLM}).
We are now able to solve the  problem of whether 
contractive projections suffice here.  This again illustrates
some limitations of the Banach algebra category when 
studying operator algebras.    
 
\begin{corollary}  \label{condex}  There exists a commutative
operator algebra $A$ with an identity of norm $1$, 
a contractive projection $P$ on $A$
whose range is a subalgebra $B$ containing $1_A$, 
and elements $a \in A, b \in B$, such that $P(b a) \neq b P(a)$.
\end{corollary}

\begin{proof}  If there existed no operator algebra with this
property, 
then it is shown in \cite{Bare} that Theorem \ref{main} 
would fail.
\end{proof}

{\bf Remark:}  
The fact that the algebra in Corollary \ref{condex}
is commutative, also appears to rule out the existence of
a `Jordan algebra variant' of Tomiyama's result for 
contractive projections, in the setting of 
nonselfadjoint operator algebras.  We thank J. Arazy for discussions
on this point in 2002; he also suggested (with a different proof) the 
following partial result (which is somewhat related to
Proposition \ref{delt}):  Namely, 
if $A$ is an operator algebra with an identity of norm $1$,
and if $P$ is a contractive projection on $A$
whose range is a subalgebra $B$ containing $1_A$,
then $P(b a) = b P(a)$ for all $a \in A, b \in \Delta(B) = B \cap B^*$.
Indeed this follows from Lemma 3.2 of \cite{BM}.

\medskip

Left multiplication by a fixed element of an operator algebra, is
an example of an {\em operator space left multiplier} (e.g.\ see \cite{BM} or
\cite[Chapter 4]{BLM} for the full definition of the latter notion).
In stark contrast to  Theorem 4.1 of \cite{BM}, which is
valid for operator spaces which have an operator space predual,
we see:

\begin{corollary} There exists an operator space $B$, which is a dual Banach space,
and a left multiplier of $B$,
which is not w*-continuous on $B$.  In fact such multipliers
may be chosen to be idempotent, or nilpotent.
\end{corollary}

\end{document}